\documentclass[letterpaper, 11 pt, conference]{ieeeconf}
\usepackage{amsmath,amssymb,euscript,psfrag,latexsym,graphicx}
\usepackage{bbm,color,amstext,wasysym,subfig,cuted,mathtools, cite}
\usepackage[normalem]{ulem}
\graphicspath{{./},{./figures/}}
\usepackage{dsfont}
\usepackage{ulem}

\newcommand{\cC}{{\cal C}}
\newcommand{\cK}{{\cal K}}

\newcommand{\mR}{{\mathbb R}}
\usepackage{algorithm}
\usepackage{algorithmic}
\usepackage[algo2e,linesnumbered,ruled]{algorithm2e}

\newcommand{\mU}{{\mathbb U}}

\newcommand{\bb}{{\mathbf b}}

\newcommand{\bs}{{\mathbf s}}
\newcommand{\bv}{{\mathbf v}}
\newcommand{\bl}{{\mathbf l}}

\newcommand{\bk}{{\mathbf k}}

\newcommand{\bK}{{\mathbf K}}

\newcommand{\bff}{{\mathbf f}}
\newcommand{\bh}{{\mathbf h}}

\newcommand{\bx}{{\mathbf x}}

\newcommand{\bPsi}{{\boldsymbol \Psi}}

\newcommand{\bX}{{\mathbf X}}
\newcommand{\bY}{{\mathbf Y}}
\newcommand{\bg}{{\mathbf g}}

\newcommand{\bL}{{\mathbf L}}

\usepackage{caption}
\newtheorem{theorem}{Theorem}

\newtheorem{assumption}{Assumption}

\IEEEoverridecommandlockouts
\IEEEoverridecommandlockouts
\begin{document}

\title{\LARGE  \bf Koopman-based Policy Iteration for Robust Optimal Control}

\author{Alexander Krolicki, Sarang Sutavani, and Umesh Vaidya% <-this % stops a space
\thanks{Financial support from NSF under grants  2031573 and NSF CPS award 1932458 is greatly acknowledged.}% <-this % stops a space
\thanks{The authors are with the Department of Mechanical Engineering, Clemson University, Clemson, SC, 29634, USA;\{akrolic,ssutava,uvaidya\}@clemson.edu}%
}

\maketitle

%%%%%%%%%%%%%%%%%%%%%%%%%%%%%%%%%%%%%%%
%%%%%%%%%%%%   Abstract   %%%%%%%%%%%%%
%%%%%%%%%%%%%%%%%%%%%%%%%%%%%%%%%%%%%%%

\begin{abstract}
Classically, the optimal control problem in the presence of an adversary is formulated as a two-player zero-sum differential game or an $H_\infty$ control problem. The solution to these problems can be obtained by solving the Hamilton-Jacobi-Issac equation (HJIE). We provide a novel Koopman-based expression of the HJIE, where the solutions can be obtained through the approximation of the Koopman operator itself. In particular, we developed a data-driven and model based policy iteration algorithm for approximating the optimal value function using a finite-dimensional approximation of the Koopman operator and generator.

\end{abstract}

\section{Introduction}

The control of complex dynamical systems in applications such as vehicle autonomy, robotics, and advanced manufacturing involves interactions of system dynamics with the environment. There are two common approaches to account for these interactions: robust control formulation using $H_\infty$ optimal control or control with adversary using two-player zero-sum differential games \cite{basar1999}. These two formulations are intimately connected and essentially reduce to solving the HJIE \cite{basar1995}. The HJIE is a nonlinear partial differential equation (PDE) that is difficult to solve for nonlinear dynamical systems. Numerical methods are often used for the approximate solution to the HJIE. In this paper, we provide a novel Koopman-based perspective to the HJIE. The Koopman operator provides a linear lifting of a nonlinear system in the space of observables or functions \cite{Lasota}. The linearity of the  Koopman operator provides a powerful tool in the development of data-driven analysis and synthesis methods for nonlinear systems \cite{brunton2016koopman,peitz2019koopman,villanueva2021towards,borggaard2009control, abraham2019active,korda2018linear, sootla2018optimal,otto2021koopman,kaiser2021data}.

Given the significance of $H_{\infty}$ optimal control and two-player zero-sum differential game problems, there is extensive literature on the computational aspects of this problem.

Many of the algorithms that we consider require an initial stabilizing controller and follow a nested loop structure \cite{Bea1998SuccessiveGA,Lewis2006Hinfty,Lewis2007HJI}.  The inner loop iterates to find the worst adversary or disturbance, given a controller.  The outer loop works to find the best response to the adversary; for more details, we refer the reader to \cite{Bea1998SuccessiveGA}. Our algorithm will utilize the same architecture. In contrast to nested loops, the algorithm in \cite{Lewis2010Synchronous} computes the adversary and control updates simultaneously. Other techniques for solving the HJIE are adapted from works developed to solve the Hamilton-Jacobi-Bellman equation (HJBE). Recent works have adapted the Kleinman algorithm \cite{KleinmanAlg} by converting the problem of solving a Riccati equation, with a sign indefinite quadratic term, into one of generating successive iterations of solutions of LQ-type Riccati equations, each with negative semidefinite quadratic terms \cite{Lanzon2007}. These results were extended for solving the HJIE with guarantees on the local convergence rate and no requirements for an initial admissible control \cite{FENG2009881}. These works and others can benefit from the Koopman-based HJIE as we provide new insights into solutions for these equations.

The main contribution of this paper is a Koopman-based derivation of the HJIE. In this paper, we provide model-based and data-driven algorithms for the approximate solution of the HJIE. Given the intimate connection between the HJIE and $H_\infty$ optimal control, our proposed Koopman-based computational framework will solve either of these problems. 

The paper is organized as follows. Preliminaries on the Koopman operator and its data-driven approximation are presented in Section II. Results involving the connection between the HJIE and the Koopman operator are presented in Section III. Model-based and data-driven algorithms for obtaining iterative solutions of the HJIE based on the approximation of the Koopman operator are presented in Section IV. Simulation results of the developed framework on nonlinear system examples are presented in Section V, followed by conclusions in Section VI. 
\section{Preliminaries and Notations}
\noindent{\bf Notation:} $\mR^n$ denotes the $n$ dimensional Euclidean space  and $\mR^n_{\geq 0}$ is the positive orthant. Given $\bX\subseteq \mR^n$ and $\bY\subseteq \mR^m$, let ${\cal L}_1(\bX,\bY), {\cal L}_\infty(\bX,\bY)$, and ${\cal C}^k(\bX,\bY)$ denote the space of all real valued integrable functions, essentially bounded functions, and space of $k$ times continuously  differentiable functions mapping from $\bX$ to $\bY$ respectively.

$\bs_t(\bx)$ denotes the solution of dynamical system $\dot \bx={\bf f}(\bx)$ starting from initial condition $\bx$. 
\subsection{Koopman Operators and Generators}
Consider a dynamical system
% \vspace{-0.1in}
\begin{align}\dot \bx={\bf f}(\bx),\;\;\bx\in \bX\subseteq \mathbb{R}^n, \label{dyn_sys}\end{align} where the vector field is assumed to be ${\bf f}(\bx)\in {\cal C}^1(\bX,\mR^n)$. The nonlinear dynamics in the state space can be lifted to infinite dimension space of functions using the Koopman operator defined as follows.\\

The Koopman operator $\mathbb{U}_t :{\cal L}_\infty(\bX)\to {\cal L}_\infty(\bX)$ for dynamical system~\eqref{dyn_sys} is defined as 
\[[\mathbb{U}_t \varphi](\bx)=\varphi(\bs_t(\bx)),\;\;\varphi\in {\cal L}_\infty. 
\]
The infinitesimal generator for the Koopman operator 
\begin{equation}
\lim_{t\to 0}\frac{\mathbb{U}_t\varphi-\varphi}{t}=\bff(\bx)\cdot \nabla \varphi(\bx)=:{\cal K}_{\bff} \varphi. \label{K_generator}
\end{equation}
Since $\mU_t$ is semi-group with generator ${\cal K}_{\bff}$ it satisfies
\begin{align}
    \frac{d}{dt}\mU_t\varphi={\cal K}_{\bff} \mU_t\varphi \label{gen}.
\end{align}

Next we describe algorithms for the finite dimensional approximation of Koopman operator.
\subsection{Data-driven Approximation of Koopman Operator}\label{section_dataapproximation}
Extended dynamic mode decomposition (EDMD) method is one of the popular algorithms for the data-driven approximation of the Koopman operator \cite{williams2015kernelbased}.  The basic idea of the algorithm can be explained as follows. For the continuous-time dynamical system (\ref{dyn_sys}), consider snapshots of time-series data from single or multiple trajectories
\begin{eqnarray}
{\mathcal X}= [\bx_1,\bx_2,\ldots,\bx_M],\;\;\;\;{\cal Y} = [\mathbf{y}_1,\mathbf{y}_2,\ldots,\mathbf{y}_M] ,\label{data}
\end{eqnarray}
where $\bx_i\in \bX$ and $\mathbf{y}_i\in \bX$. The pair of data sets are assumed to be two consecutive snapshots i.e., $\mathbf{y}_i=\bs_{\Delta t}(\bx_i)$, where $\bs_{\Delta t}$ is solution of (\ref{dyn_sys}) with $\Delta t$ the discretization time-step. Let ${\bPsi}=[\psi_1,\ldots,\psi_N]^\top$ be the choice of basis functions.
The EDMD algorithm provides a finite-dimensional approximation of the Koopman operator as the solution of the following least square problem $\min\limits_{\bf K}\parallel {\bf G}{\bf K}-{\bf A}\parallel_F^2$.
With ${\bf K},{\bf G},{\bf A}\in\mathbb{R}^{N\times N}$, $\|\cdot\|_F$ is the Frobenius norm. The matrices are computed from data ${\bf G}=\frac{1}{M}\sum_{m=1}^M \bPsi({\bx}_m)\bPsi({\bx}_m)^\top$ and ${\bf A}=\frac{1}{M}\sum_{m=1}^M \bPsi({\bx}_m) \bPsi({\mathbf y}_m)^\top$.
The solution to the least square problem is
\begin{eqnarray}
\mU_{\Delta t}\approx{\bf K}_{edmd}=\bf{G}^\dagger \bf{A}\label{edmd_formula}.
\end{eqnarray}
Where $\dagger$ stands for pseudo-inverse. The convergence of EDMD towards the true Koopman operator as the number of data points and basis functions go to infinity are provided in \cite{korda2018convergence,KLUS2020132416}. The EDMD-based approximation of the Koopman operator can be used to approximate the Koopman generator as follows,
\begin{align}\label{Gene_def}
    {\cal K}_{\bff}\approx \frac{\bK-I}{\Delta t}=: \bL.
\end{align}

\section{Koopman and Hamilton-Jacobi Equations}
There are different ways in which control problems in the presence of an adversary can be formulated.
In particular, the adversary can be viewed as a passive disturbance as in the $H_\infty$ problem, or it can be viewed as an active player as in the two-player zero-sum differential game problem. In this paper, we use the latter approach. Consider an affine in control and disturbance system of the form
\begin{align}
    \dot \bx=\bff (\bx)+\bg(\bx)u(t)+\bh(\bx) \omega(t), \label{System_dynamics}
\end{align}
where $\bx\in \mR^n$, $u\in \mR$, and $\omega\in \mR$ are the state, control input, and adversarial input respectively. For the simplicity of presentation, we restrict the discussion to the case of scalar control and adversary inputs. We make the following assumptions on the vector fields. 
\begin{assumption}\label{assume_dynamics}
We assume that the vector fields $\bff,\bg$, and $\bh$ are $\cC^1(\mR^n)$ and that $\bff(0)=0$ so that the origin is the equilibrium point of the system in the absence of control and adversarial inputs. 
\end{assumption}

In the two-player zero-sum differential game, we proceed with the following performance index. 
\begin{equation}
\resizebox{.4\textwidth}{!}{$
\begin{split}
    J(\bx,u,w)=\int_0^\infty \left(q(\bx)+r u^2(t)-\gamma^2 \omega^2(t)\right) dt, \label{perf_cost}
\end{split}$
}
\end{equation}
where $q(\bx)\geq 0$ is assumed to be a state cost with $q(\bx)\in \cC^1(\mR)$ and $\gamma,r$ are positive constants. The two-player zero-sum differential game is defined as the following min-max optimization problem
\begin{align}
    \min_{u[0,\infty)}\max_{w[0,\infty)} J(\bx,u,w),
\end{align}
where the objective of the control and adversarial inputs is to minimize and maximize the performance, respectively. 
% This is termed as the two-player zero-sum game because the control input and the adversary input compete to drive the performance index in the opposite directions. 
The game yields a unique solution $V^{\star}$ (known as the value of the game) if a saddle point $(u^{\star},w^{\star})$ exists that satisfies the following condition, $V^\star(\bx)=\min _{u} \max _{w} J(\bx, u, w)=\max _{w} \min _{u} J(\bx, u, w)$.
This saddle point follows the so called Nash equilibrium condition, namely
$J(\bx, u^{\star}, w) \leq J\left(\bx, u^{\star}, w^{\star}\right) \leq J\left(\bx, u, w^{\star}\right)$.
It is known that the optimal value function $V^\star(\bx)$ is obtained from the solution of the HJIE 
\begin{align}
0= q+  \bff(\bx)\cdot \nabla V-
\frac{r^{-1}}{4} \nabla V^{\top} \bg(x)  \bg^{\top}(\bx) \nabla V \nonumber\\
+\frac{1}{4 \gamma^{2}} \nabla V^{\top} \bh(x)  \bh^{\top}(\bx) \nabla V, \quad V(0)=0 \label{HJI_orig}
\end{align}
and, the optimal control and adversarial input obtained from $V^\star$ is as follows:
\begin{equation}
u(\bx)= -\frac{1}{2} r^{-1} \bg\cdot \nabla V^\star, \ w(\bx)=\frac{1}{2 \gamma^{2}} \bh\cdot \nabla V^\star.
\end{equation}
The minimum positive semi-definite (PSD) solution gives the value of the game or the Nash value. It has been shown that a unique minimum PSD solution exists for all $\gamma \geq \bar \gamma$, where the $\bar \gamma$ corresponds to the solution of the $H_{\infty}$ control problem. For more details on this connection between the differential games and the $H_\infty$ problem, refer to \cite{basar1995}. 

In the following, we derive the HJIE using the Koopman theory. This result is novel and will lay the foundation for the further development of Koopman-based computational frameworks. Consider the min-max performance index subject to the constraints of system dynamics as follows
\begin{align}
    &\min_u \max_w J(\bx,u,w)\nonumber\\
    &{\rm s.t.}\;\; \dot \bx=\bff (\bx)+\bg(\bx)u(t)+\bh(\bx) \omega(t), \label{minmax_prob}
\end{align}
where $J$ is defined in (\ref{perf_cost}). We make following assumption on the min-max problem (\ref{minmax_prob}). 
\begin{assumption}\label{assume_optimization} We assume that the solution to the min-max problem (\ref{minmax_prob}) for the control and adversarial inputs is feedback in nature i.e., $u=k(\bx)\in \cC^1(\mR^n)$ and $w=\ell(\bx)\in{\cC}^1(\mR^n)$. Furthermore, $\gamma$ is assumed to be larger than $\bar \gamma$ so that the optimal value function, $V^\star(\bx)\geq 0$, is finite for any finite value of $\bx\in \mR^n$. 
\end{assumption}
Following Assumption \ref{assume_optimization}, we can write  
\begin{gather}
   V^\star(\bx)= \min_k \max_\ell J(\bx,k,\ell)\nonumber\\
    {\rm s.t.}\;\; \dot \bx=\bff (\bx)+\bg(\bx)k(\bx)+\bh(\bx)\ell(\bx). \label{minmax_prob2}
\end{gather}
The above assumption could be restrictive, and in fact, the HJIE for solving the min-max problem is derived under a less restrictive assumption \cite{basar1995}. For example, the assumption rules out the possibility for the existence of viscosity based solution of the HJIE \cite{crandall1992user}. Essentially, the viscosity-based solutions allow for a continuous function to be defined as a unique solution of the HJIE.

\begin{theorem}
For the min-max optimization problem (\ref{minmax_prob}) satisfying Assumption \ref{assume_optimization} and system dynamics Assumption \ref{assume_dynamics}. The optimal cost function $V^\star$ can be obtained as the solution of following equation
\begin{equation}\label{hji} 
\begin{split}
&{\cal K}_{\bff+\bg k+\bh\ell} V^\star=-q-r(k^\star)^2 +\gamma^2(\ell^\star)^2  \\
&k^\star(\bx)= -\frac{1}{2} r^{-1} {\cal K}_{\bg}V^\star , \ \ell^\star(\bx)=\frac{1}{2 \gamma^{2}} {\cal K}_{\bh} V^\star,
\end{split}
\end{equation}
where  ${\cal K}_{\bff+\bg k+\bh\ell}$, ${\cal K}_{\bg }$, and ${\cal K}_{\bh}$ are the Koopman generators where the subscript denotes the vector field.  
\end{theorem}
\begin{proof}Following Assumption \ref{assume_optimization}, we write the feedback system as 
\[\dot \bx=\bff(\bx)+\bg(\bx)k(\bx)+\bh(\bx)\ell(\bx).\]
Let $\mU_t^c$ be the Koopman operator for this feedback system. Using the definition of Koopman operator, we can write performance measure $J(\bx,k,\ell)$ as 
\begin{align}V(\bx)=\int_0^\infty [\mU_t^c\varphi](\bx) dt, \label{integral}
\end{align}
where $\varphi(\bx):=(q+r k^2-\gamma^2 \ell)(\bx)$. We next claim that $V(\bx)$ satisfies
\begin{align}{\cal K}_{\bff+\bg k+\bh\ell}V=-q-rk^2 +\gamma^2\ell^2. \label{ee}
\end{align}
Substituting (\ref{integral}) in the LHS of (\ref{ee}), we obtain
\begin{align}
\int_0^\infty {\cal K}_{\bff+\bg k+\bh\ell}\mU_t^c \varphi dt= \int_0^\infty  \frac{d}{dt}\mU_t \varphi dt \nonumber\\
=\mU_t\varphi|_{t=0}^\infty=\lim_{t\to \infty}\mU_t \varphi-\varphi, 
\end{align}
where we have used the infinitesimal generator property of the Koopman semi-group (\ref{gen}). Then, we can show that $\lim_{t\to \infty}[\mU_t \varphi](\bx)=0$. To prove this we use the Assumption \ref{assume_optimization} that $V(\bx)$ is finite for any finite $\bx$. Furthermore, $[\mU_t \varphi]$ is uniformly continuous w.r.t. time which follows from the definition of the Koopman operator semi-groups and the fact that the solution of the closed loop system is uniformly continuous w.r.t. time.
Hence, we can apply the Barbalat Lemma which states that for any function $G(t)\in {\cal C}^1$, and $\lim_{t\to \infty} G(t)=\alpha$. If $G'(t)$ is uniformly continuous, then $\lim_{t\to \infty} G'(t)=0$. Applying the Barbalat Lemma with $G'(t)=[\mU_t \varphi](\bx)$ for a fixed $\bx$, we obtain $\lim_{t\to \infty}[\mU_t \varphi](\bx)=0$. Hence, we prove the claim (\ref{ee}). The optimal control and adversarial inputs are obtained as the critical or extremum point of (\ref{ee}). This extremum is obtained by differentiating  (\ref{ee}) w.r.t. $k$ and $\ell$ leading to
\begin{eqnarray}
k(\bx)= -\frac{1}{2} r^{-1} {\cal K}_{\bg} V,\;\;\;\;
\ell(\bx)=\frac{1}{2 \gamma^{2}} {\cal K}_{\bh} V. \label{control_form}
\end{eqnarray}
Substituting (\ref{control_form}) in (\ref{ee}) we obtain the desired HJIE (\ref{hji}) for optimal value function $V^\star$. 
\end{proof}

Note that the Koopman-based perspective developed here can be easily extended to the HJBE and the $H_\infty$ control problem. The HJBE appears in the optimal control problem in absence of the adversary and can be obtained from (\ref{HJI_orig}) when $\gamma\to \infty$.

\section{Computational Methods}
The complexity associated with the nonlinear nature of the PDE is overcome by developing an iterative algorithm for solving the HJIE. In the following, we present the Koopman Policy Iteration (KPI) algorithm.  
\subsection{Koopman policy iteration approximation of HJIE}\label{section_modelbased}
In the approximation of the HJIE, we can assume that we have access to the system vector fields in the form of $\bff,\bg$, and $\bh$. If we do not, then EDMD based methods can be employed to approximate these vector fields. For a given control and adversarial input at the iteration step $(i,j)$ i.e., $k^{i}$ and $\ell^{(i,j)}$, the time-series data from system dynamics
$\dot \bx=\bff+\bg k^i+\bh\ell^{(i,j)}=: \bff_c^{(i,j)}$ is used to construct the approximation of the closed loop Koopman generator $\cK_{\bff_c^{(i,j)}}\approx \bL_{(i,j)}^{\bff_c}$ for a given choice of basis functions $\bPsi=(\psi_1,\ldots, \psi_N)^\top$.
Similarly, the right hand side of (\ref{ee}) can be approximated as follows:
\begin{align}
    -q(\bx)-r (k^i)^2(\bx)+\gamma^2(\ell^{(i,j)})^2(\bx)\approx \bb_{(i,j)}^\top \bPsi(\bx),\nonumber
\end{align}
 where the coefficient vector $\bb_{(i,j)}$ is obtained as the solution of following least square problem.
Let   $\{\bx_t\}_{t=0}^M$ be the time series data collected by sampling uniformly over the state space. Construct the matrices for the lifted states $\bar{\bPsi} = [ \bPsi(\bx_0), ..., \bPsi(\bx_M)]^\top$, state cost $\bar{\mathbf{q}} = [q(\bx_0), ..., q(\bx_M)]^\top$, control cost $\bar{\bk} = [(k^i)^2(\bx_0), ..., (k^i)^2(\bx_M)]^\top$, and adversary cost $\bar{\bl} = [(\ell^{(i,j)})^2(\bx_0), ..., (\ell^{(i,j)})^2(\bx_M)]^\top$. 
Where the coefficient vector $\bb_{(i,j)}$ is obtained as the solution of following least square problem,
\begin{equation}\label{ls_cost}
\min\limits_{\bf \bb_{(i,j)}}\parallel \bar{\bPsi}(\bx)\bb_{(i,j)}-(-\bar{\mathbf{q}}-r \bar{\bk}^i+\gamma^2\bar{\bl}^{(i,j)})\parallel_2^2,
\end{equation}
which admits the following analytical solution.
\begin{equation}\label{ls_cost_sol}
\bb_{(i,j)}= \bar{\bPsi}^\dagger(-\bar{\mathbf{q}}-r \bar{\bk}^i+\gamma^2\bar{\bl}^{(i,j)}).
\end{equation}
The value function is approximated as $V^{(i,j)}(\bx)=\bv_{(i,j)}^\top \bPsi(\bx)$, then the finite dimensional approximation
\begin{figure}[H]
\centering
\includegraphics[width=0.45\textwidth]{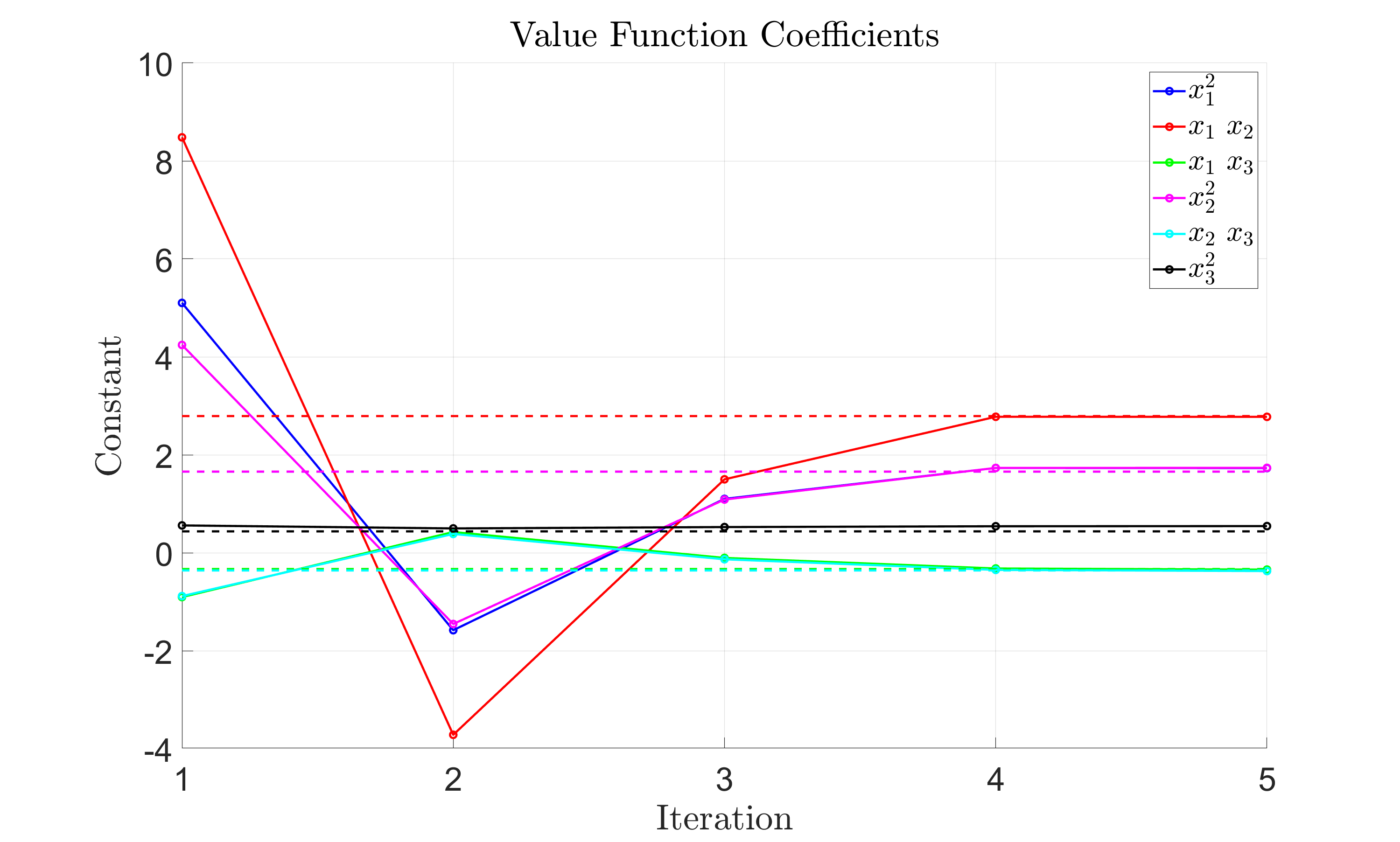}
\caption{Example 1: Analytical (dashed) and KPI (solid).}
\label{fig:3D MB value}
\includegraphics[width=0.45\textwidth]{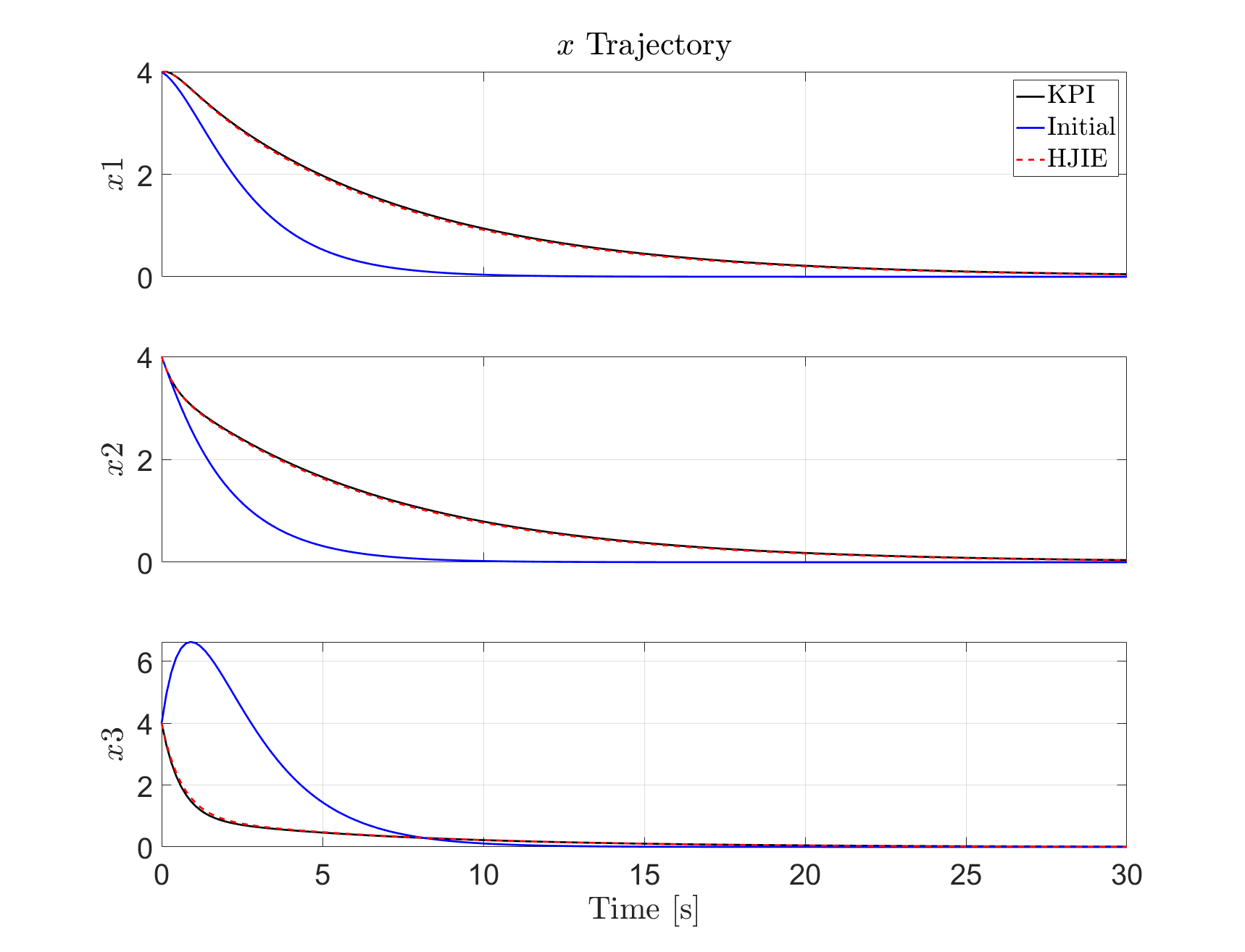}
\caption{Example 1: State trajectories.}
\label{fig:3D MB state}
\includegraphics[width=0.45\textwidth]{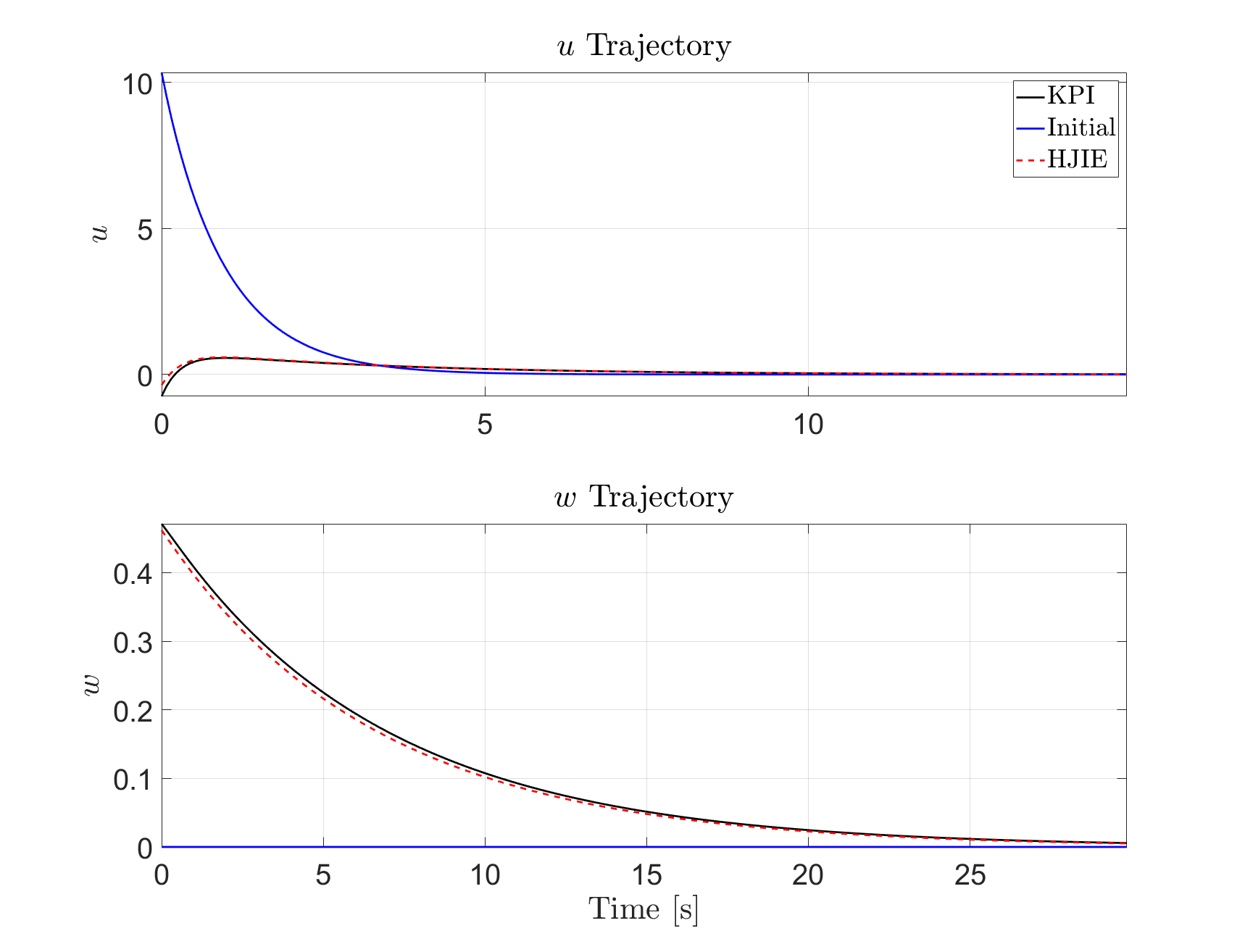}
\caption{Example 1: Control and adversary trajectories.}
\label{fig:3D MB ca}
\end{figure}
of (\ref{hji}) can be written as 
\begin{align}
    \bL^{\bff_c}_{(i,j)}\bv_{(i,j)}=\bb_{(i,j)}, \label{aa1}\\
    k^{(i+1)}(\bx)=-\frac{1}{2}r^{-1}\bg\cdot \nabla (\bv^{\top}_{(i,\infty)}\bPsi(\bx)), \\
     \ell^{(i,j+1)}(\bx)=\frac{1}{2\gamma^2}\bh\cdot \nabla (\bv^{\top}_{(i,j)}\bPsi(\bx)).
\end{align}

\section{Simulation Results}
All the simulation results assume knowledge of the vector fields $\bff,$ $\bg,$ and $\bh$. All the computations are performed in MATLAB with a i9-10900KF CPU.
\subsection {Example 1: F16 Aircraft}
%%%%%%%%%%%%%%%%%%%% 3D System %%%%%%%%%%%%%%%%%%%%
The HJIE for linear systems can be solved analytically, so we validate our KPI algorithms solution. The system dynamics of a F16 aircraft is given as a linear continuous time model $\dot{x} = Ax+Bu+H\omega$,
\begin{equation} \label{3D sys}
\resizebox{.425\textwidth}{!}{
$\dot{x}=\begin{bmatrix} 
        -1.01887 & 0.90506 & -0.00215\\
        0.82225 & -1.07741 & -0.17555\\
        0 & 0 & -1 
        \end{bmatrix} x +
        \begin{bmatrix} 
        0 \\
        0\\
        1 
        \end{bmatrix} u + \begin{bmatrix} 
        1 \\
        0\\
        0 
        \end{bmatrix} \omega$.
}
\end{equation}
The basis functions used in this problem are $\Psi(x)=[x_1^2, x_1 x_2, x_1 x_3, x_2^2,  x_2 x_3, x_3^2]$. The time step between snapshots is $\Delta T = 0.15s$ with a total of 10 initial conditions sampled uniformly from $x_{ic} \in [-5,5]$. We choose $\gamma = 5$ and $r=1$. Let the initial admissible control using pole placement method with poles $p=[-2, -1, -0.5]$.
The game algebraic Ricatti equation (GARE) for $H_\infty$ control problem takes the form, 
\begin{equation} \label{hcare}
\resizebox{.425\textwidth}{!}{$
A^\top P+PA-P[H,B]\begin{bmatrix} -\gamma^2I & 0 \\
                                    0       &   I      
                                    \end{bmatrix}
                                    \begin{bmatrix}
                                            H^T \\
                                            B^T
                                    \end{bmatrix}P +
                                    C^\top C = 0. \nonumber
$}
\end{equation}
The analytic value function, using the GARE solution $P$, is $V^{\star}=x^\top Px$. In this example, the $P$ matrix is, 
\begin{equation} \label{3D P}
\resizebox{.25\textwidth}{!}{$
\begin{split}
        \begin{bmatrix} 
        1.657 & 1.395 & -0.166\\
      1.395 & 1.657 & -0.180\\
        -0.166 & -0.180 & 0.437 
        \end{bmatrix}.
\end{split}$
}
\end{equation}
Our proposed KPI algorithm converges to coefficient vector
$\bv  = [p_{11}, 2p_{12}, 2p_{13}, p_{22}, 2p_{23}, p_{33}]^\top = [1.7293, 2.7756, -0.343, 1.733, -0.373, 0.544]^\top$, which closely matches with the analytical solution obtained using GARE. The KPI value function coefficients converge to the HJIE solutions in 4 iterations, see Fig. \ref{fig:3D MB value}. The closed loop state trajectories from the initial, KPI, and $H_\infty$ control is shown Fig. \ref{fig:3D MB state}, and \ref{fig:3D MB ca}.
The total algorithm runtime was 0.243s. 
\begin{figure}[t]
\includegraphics[width=0.5\textwidth]{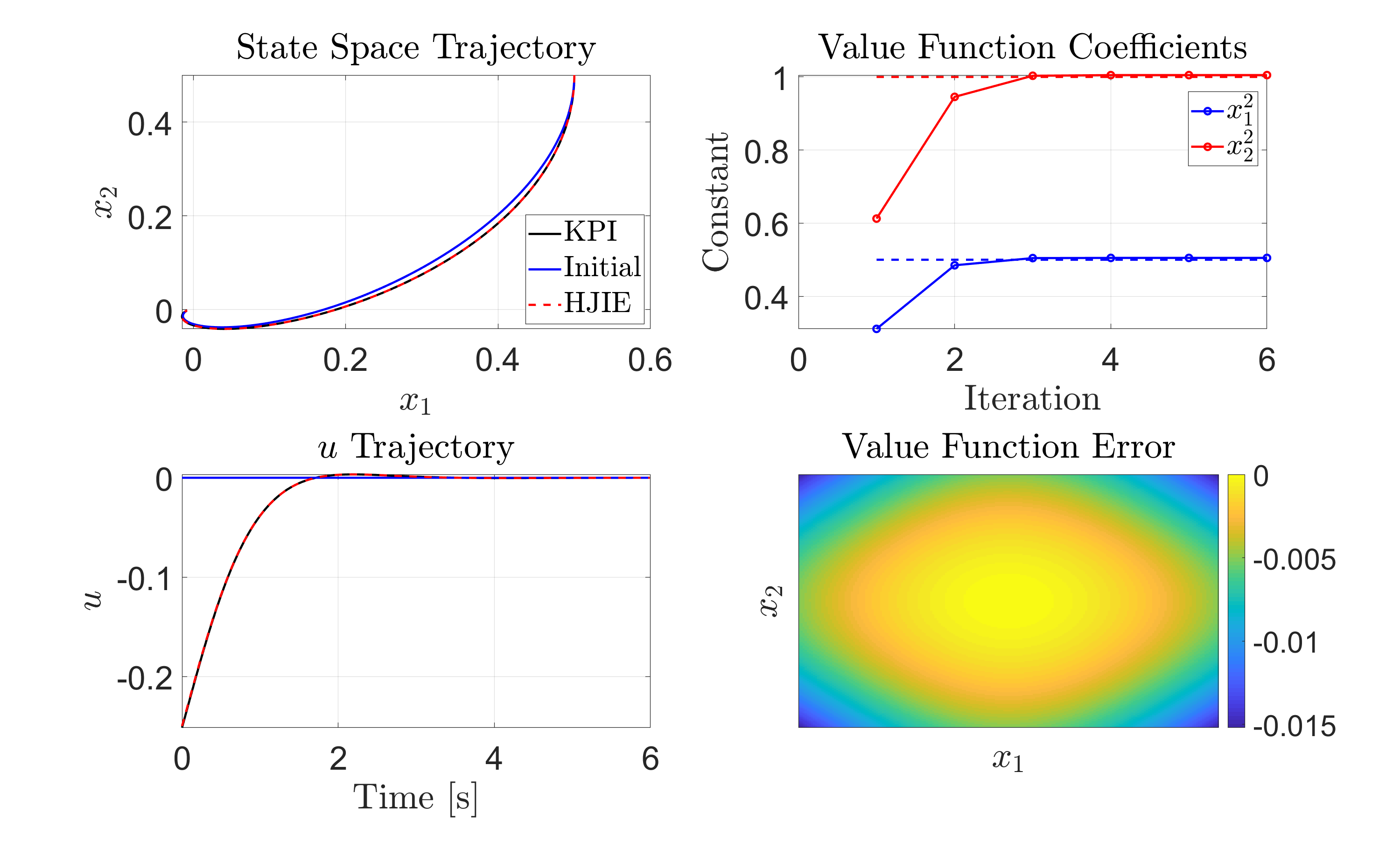}
\caption{Example 2: HJBE simulation results. Top left: state space trajectory, bottom right: control input trajectory, top right: value function coefficients (dashed lines are analytical values), and  bottom right: value function error.}
\label{fig:2D MB state quad}
\vspace{5.00mm} 
\includegraphics[width=0.5\textwidth]{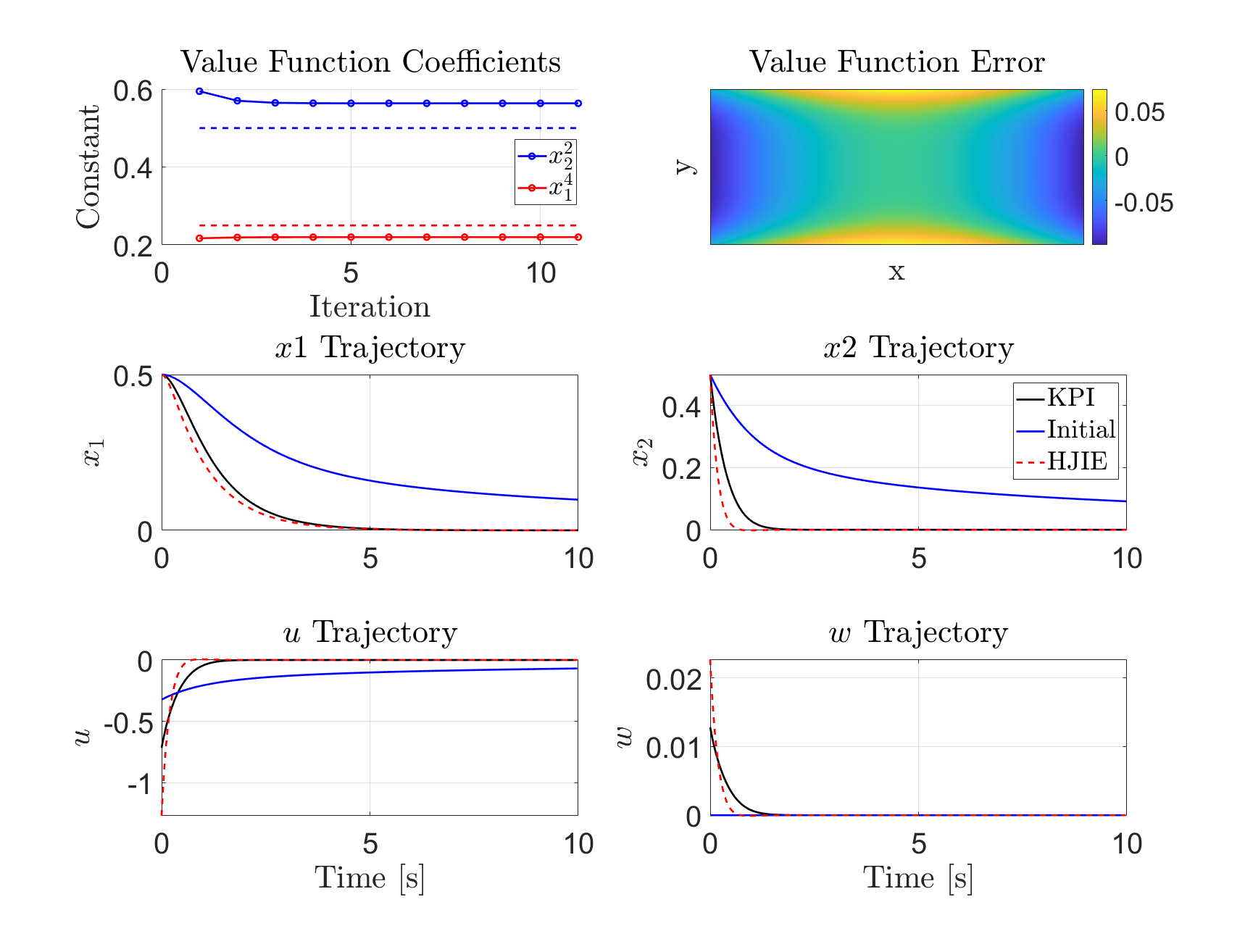}
\caption{Example 3: HJIE simulation results. Top row: (left) value function coefficients and (right) error, middle row: state response over time, and bottow row: (left) control and (right) adversary input values over time.}
\label{fig:2D state quad}
\end{figure}
\subsection {Example 2: 2D Nonlinear System HJBE}
%%%%%%%%%%%%%%%%%%%% 2D System %%%%%%%%%%%%%%%%%%%% 
The continuous time nonlinear dynamical system is,
\begin{equation} \label{3D sys}
\resizebox{.425\textwidth}{!}{
$\dot{x}=\begin{bmatrix} 
        -x_1 + x_2 \\
        -0.5(x_1+x_2)+0.5(x_1^2x_2) 
        \end{bmatrix} + 
        \begin{bmatrix} 
        0 \\
        x_1
        \end{bmatrix} u + \begin{bmatrix} 
        0.1 \\
        0.5x_2
        \end{bmatrix} \omega$.
}
\end{equation}
We choose polynomial basis functions of the form $\bPsi(\bx)=[ x_1^2 x_2, x_1^2, x_1 x_2, x_2^2, x_1^2 x_2^2, x_1^4,x_2^4]$. The initial admissible control is chosen to be $u_{initial}=0$, as the system has a stable equilibrium at the origin. In absence of $\omega$ or when $\gamma=\infty$, the analytical solution of HJBE for the optimal value function is $V^{\star}=0.5x_1^2+x_2^2$, hence $\bv^\star(x_1^2) = 0.5$ and $\bv^\star(x_2^2) = 1$. The HJBE optimal control signal is $u^\star(x)=-x_1x_2$. The analytical solutions can be found in \cite{Vamvoudakis2009OnlineAC}. In this case $r=1$ and $\Delta T = 0.01s$ are used with a total of 50 initial conditions sampled from $x_{ic} \in [-1,1]$. The KPI value function for this case converges to coefficient vector $\bv = [0, 0.505, 0, 1.00, 0.005, 0, 0]^\top$, which matches the analytical solution. The solution converges to the optimal HJBE control in 3 iterations as shown in Fig. \ref{fig:2D MB state quad} with a total run-time of 0.396s. 
In the second case, $\omega$ is introduced and $\gamma=5$ is chosen. The optimal value function is not know. The KPI algorithm converges to a value function with coefficients $
\bv=[0, 0.504, -0.002, 1.005, 0.003, 0.001, 0.003]^\top$.
\subsection{Example 3: 2D Nonlinear System HJIE}
%%%%%%%%%%%%%%%%%%%% 2D System %%%%%%%%%%%%%%%%%%%% 
The continuous time nonlinear dynamical system is,
\begin{equation} \label{3D sys}
\resizebox{.425\textwidth}{!}{$
\begin{split}
\dot{x}=\begin{bmatrix} 
        -x_1 + x_2 \\
        -x_1^3-x_2^3+\frac{x_2(\cos(2x_1)+2)^2}{4}-\frac{x_2(\sin(4x_1)+2)^2}{4\gamma^2}
        \end{bmatrix} + \\
        \begin{bmatrix} 
        0 \\
        \cos(2x_1)+2
        \end{bmatrix} u + \begin{bmatrix} 
        0 \\
        sin(4x_1)+4\\
        \end{bmatrix} \omega.
\end{split}$
}
\end{equation}
We use the following basis $\Psi(x)=[x_1, x_2, x_1^{2}x_2, x_1^3, x_2^2, x_2^3, x_1^4]$. In this example, $r=1$, $\gamma=8$ and, $\Delta T = 0.025s$ with a total of 100 initial conditions sampled from $x_{ic} \in [-1.25,1.25]$. $u = \frac{(\frac{x_2(\cos(2x_1)+2)^2}{4}-\frac{x_2(\sin(4x_1)+2)^2}{4\gamma^2})}{(\cos(2x_1)+2)}$ is used as the initial controller. The solution to the HJIE for the optimal value function is $V^\star=\frac{1}{4}x_1^4+\frac{1}{2}x_2^2$. The optimal control is $u^\star=-(\cos(2x_1)+2)x_2$ and the optimal adversarial disturbance $\omega^\star=\frac{1}{\gamma^2}(\sin(4x_1)+2)x_2$. The analytical solution can be found in \cite{Lewis2010Synchronous}. The value function coefficients converge to $
\bv = [0.033, 0.024, -0.003, -0.023,  0.588, -0.008, 0.301]^\top$.
The solutions converges near to the HJIE in a few iterations. A more appropriate choice of basis function can help improve this approximation. The total run-time was 0.826s.
\section{Conclusions}
We presented a Koopman-based policy iteration (KPI) algorithm for solving the HJIE. This iterative algorithm can be implemented in the model-based and data-driven setting and relies on approximation of the Koopman operator. Future research efforts will focus on understanding the role of the Koopman spectrum in the data-driven approximation of the HJIE.

\bibliographystyle{IEEEtran}
\bibliography{ref1,ref}
\end{document}